\documentclass[a4paper,11pt]{amsart}

\usepackage{graphicx}
\usepackage{mathptmx}
\usepackage{amsmath}
\usepackage{amssymb}
\usepackage{enumitem}
\usepackage{xcolor}
\usepackage{pgfplots}

\newmuskip\pFqmuskip

\newcommand*\pFq[6][8]{%
  \begingroup 
  \pFqmuskip=#1mu\relax
  \mathcode`=\string"8000
  \begingroup\lccode`\~=`\,
  \lowercase{\endgroup\let~}\pFqcomma
  F^{#2}_{#3}{\left(\genfrac..{0pt}{}{#4}{#5}\bigg|#6\right)}%
  \endgroup
}
\newcommand{\pFqcomma}{\mskip\pFqmuskip}

\newtheorem{theorem}{Theorem}[section]
\newtheorem{lemma}[theorem]{Lemma}
\newtheorem{corollary}[theorem]{Corollary}

\newtheorem{remark}[theorem]{Remark}

\begin{document}

\title[]{Several expressions for degenerate harmonic numbers and some related numbers}

\author{Taekyun  Kim}
\address{Department of Mathematics, Kwangwoon University, Seoul 139-701, Republic of Korea}
\email{tkkim@kw.ac.kr}
\author{Dae San  Kim}
\address{Department of Mathematics, Sogang University, Seoul 121-742, Republic of Korea}
\email{dskim@sogang.ac.kr}

\author{Kyo-Shin Hwang}
\address{Graduate School of Education, Yeungnam University, Gyeongsan 38541, Republic of Korea}
\email{kshwang@yu.ac.kr}

\subjclass[2010]{11B73; 11B83}
\keywords{degenerate harmonic numbers}

\begin{abstract}
Many authors have recently studied the degenerate harmonic numbers. This paper makes two main contributions. First, we derive several explicit expressions for these numbers, which are a degenerate version of the ordinary harmonic numbers. We also examine the degenerate harmonic numbers of order $m$ and find an expression for them. Second, we investigate some related numbers that are closely connected to the degenerate harmonic numbers of order $m$, which reduce to the degenerate harmonic numbers when $m$=1.
\end{abstract}

\maketitle

\markboth{\centerline{\scriptsize Several expressions for degenerate harmonic numbers and some related numbers}}
{\centerline{\scriptsize Taekyun Kim, Dae San Kim and Kyo-Shin Hwang}}

\section{Introduction}
Carlitz initiated the study of degenerate Bernoulli and Euler numbers and polynomials as degenerate versions of their classical counterparts (see [5]). Recently, this pioneering work has regained the interest of some mathematicians, leading to investigations into various degenerate versions of special numbers and polynomials. This research has yielded interesting arithmetical and combinatorial results, as well as applications in other disciplines (see [3,9-16] and the references therein). For example, Kim and Kim investigated the degenerate harmonic and hyperharmonic numbers, which are degenerate versions of the classical harmonic and hyperharmonic numbers (see [11-14,16]). \par
This paper has two main goals. First, we derive several explicit expressions for the degenerate harmonic numbers, $H_{n,\lambda}$ (see \eqref{5}), as well as for the degenerate harmonic numbers of order $m$, $H_{n,\lambda}^{(m)}$ (see \eqref{33}).
We also introduce a new sequence of numbers, $K_{n,\lambda}^{(m)}$ (see \eqref{28}), which are closely related to $H_{n,\lambda}^{(m)}$ and reduce to $K_{n,\lambda}^{(1)}=H_{n,\lambda}$, for $m=1$. The motivation for introducing these new numbers is straightforward. The generating function for  $H_{n,\lambda}^{(m)}$ is obtained by replacing $\mathrm{Li}_{1,\lambda}(t)$ with $\mathrm{Li}_{m,\lambda}(t)$ in $\frac{1}{1-t}\mathrm{Li}_{1,\lambda}(t)=\sum_{n=1}^{\infty}H_{n,\lambda}t^{n}$ (see \eqref{10-1}, \eqref{33}), where
$\mathrm{Li}_{k,\lambda}(t)$ is the degenerate polylogarithm (see \eqref{10}).
Similarly, the generating function for $K_{n,\lambda}^{(m)}$ is obtained by replacing
$\mathrm{Li}_{1,-\lambda}\big(-\frac{t}{1-t}\big)$ with $\mathrm{Li}_{m,-\lambda}\big(-\frac{t}{1-t}\big)$ in $-\frac{1}{1-t}\mathrm{Li}_{1,-\lambda}\big(-\frac{t}{1-t}\big)=\sum_{n=1}^{\infty}H_{n,\lambda}t^{n}$ (see \eqref{10-1}, \eqref{28}).  \par

This paper is organized into two parts. Section 1 begins by recalling foundational concepts: harmonic numbers, harmonic numbers of order $\alpha$, unsigned Lah numbers, and derangement numbers. We then review degenerate exponentials and their compositional inverses, the degenerate logarithms, as well as degenerate polylogarithms. The section concludes with a discussion of degenerate harmonic numbers, degenerate Stirling numbers of the first kind, and degenerate derangement numbers. The main results of the paper are presented in Section 2. Here, we derive new expressions for the degenerate harmonic numbers (Theorems 2.1, 2.10, 2.12) and the degenerate harmonic numbers of order $m$ (Theorem 2.9). We also introduce the numbers $K_{n,\lambda}^{(m)}$ and deduce explicit formulas for them (Theorems 2.6, 2.7, 2.8). As general references for this paper, the reader may refer to [6-8,17-19].
\par
\vskip0.1in
It is well known that the harmonic numbers are defined by
\begin{equation}\label{1}
H_0=1, \ H_n = 1+\frac12+\frac13+\cdots +\frac1n,\ (n\in \mathbb{N}),\quad (\rm see\ [4,6-8,12]).
\end{equation}
From \eqref{1}, we note that
\begin{equation}\label{2}
 \frac1{1-t}\log\left(\frac1{1-t}\right)=\sum_{n=1}^\infty H_n \ t^n ,\quad(\rm see\ [4,6-8,12]).
\end{equation} \par
For any nonzero $\lambda\in \mathbb{R}$, the degenerate exponentials are defined by
\begin{equation}\label{3}
 e_\lambda^x (t) = \sum_{k=0}^\infty (x)_{k,\lambda}\ \frac{t^k}{k!}, \quad e_\lambda (t) = e_\lambda^{1} (t) , \quad(\rm see\ [9-16]),
\end{equation}
where
\begin{equation*}
(x)_{0,\lambda}=1,\ (x)_{n,\lambda}=x(x-\lambda)(x-2\lambda)\cdots\big(x-(n-1)\lambda\big),\ (n\ge 1). 	
\end{equation*}
The compositional inverse of $e_\lambda(t)$ is denoted by $\log_{\lambda}(t)$ and called the degenerate logarithm. Then we have
\begin{equation}\label{4}
\log_\lambda(1+t) =\sum_{n=1}^\infty \lambda^{n-1} \ (1)_{n,1/\lambda}\ \frac{t^n}{n!}=\sum_{n=1}^{\infty}\binom{\lambda-1}{n-1}\frac{t^{n}}{n},\quad(\rm see\ [10]).
\end{equation}
From \eqref{4}, we note that $e_\lambda(\log_\lambda(1+t))=\log_\lambda(e_\lambda(1+t))= 1+t$. \par
Recently, Kim-Kim introduced the degenerate harmonic numbers given by (see [11-14,16])
\begin{equation}\label{5}
H_{0,\lambda}=0,\ H_{n,\lambda}=\frac1\lambda\sum_{k=1}^n \binom{\lambda}{k}(-1)^{k-1}=\sum_{k=1}^{n}\binom{\lambda-1}{k-1}\frac{(-1)^{k-1}}{k}, \ (n \in \mathbb{N}).
\end{equation}
Note that
$$\lim_{\lambda\to 0} H_{n,\lambda}=H_n.$$
From \eqref{5}, we have
 \begin{equation}\label{6}
 \frac1{1-t}\left(-\log_{\lambda}(1-t)\right) =\frac1{1-t}\log_{-\lambda}\left(\frac1{1-t}\right)
 =\sum_{n=1}^\infty H_{n,\lambda}\  t^n, \quad(\rm see\ [11]).
\end{equation} \par
The unsigned Lah number $L(n,k)$ counts the number of ways a set of $n$ elements can be
partitioned into $k$ nonempty linearly ordered subsets. Explicitly, the unsigned Lah number $L(n,k)$ is given by the formula (see [6,8,18]):
\begin{equation}\label{7}
L(n,k)=\binom{n-1}{k-1}\frac{n!}{k!}, \quad (n\ge k\ge 1).
\end{equation}
From \eqref{7}, we note that
\begin{equation}\label{8}
\frac1{k!}\left(\frac{t}{1-t}\right)^k =\sum_{n=k}^\infty L(n,k)\ \frac{t^n}{n!}, \quad (k\ge 0).
\end{equation} \par
For $k \in \mathbb{Z}$, the polylogarithm is defined by
\begin{equation}\label{9}
\mathrm{Li}_k (t) =\sum_{n=1}^\infty \frac{t^n}{n^k}, \quad (|t|<1), \quad(\rm see\ [10]).
\end{equation}
Note that $\mathrm{Li}_1 (t)= -\log (1-t)$. For $k \in \mathbb{Z}$, Kim-Kim introduced the degenerate polylogarithm given by
\begin{equation}\label{10}
\mathrm{Li}_{k,\lambda}(t) =\sum_{n=1}^\infty \frac{(-\lambda)^{n-1} \ (1)_{n,1/\lambda}}{(n-1)!\ n^k } \ t^n , \quad (|t|<1),\quad(\rm see\ [10]).
\end{equation}
Note that
\begin{align}
&\lim_{\lambda\to 0}\mathrm{Li}_{k,\lambda}(t)=\mathrm{Li}_{k}(t),\quad \mathrm{Li}_{1,\lambda}(t) = -\log_\lambda (1-t), \label{10-1} \\
&\frac{1}{1-t}\mathrm{Li}_{1,\lambda}(t)=-\frac{1}{1-t}\mathrm{Li}_{1,-\lambda}\big(-\frac{t}{1-t}\big)=\sum_{n=1}^{\infty}H_{n,\lambda}t^{n}. \nonumber
\end{align}
The harmonic numbers of order $\alpha$ are defined by
\begin{equation}\label{11}
H_0^{(\alpha)}=0,\quad H_{n}^{(\alpha)}=\sum_{k=1}^n\frac1{k^\alpha},\quad(\rm see\ [12]).
\end{equation} \par
The degenerate Stirling numbers of the first kind are defined by
\begin{equation}\label{12}
(x)_{n}=\sum_{k=0}^{n}S_{1,\lambda}(n,k)\ (x)_{k,\lambda},\quad (n\ge 0),\quad(\rm see\ [2,10,11,13]).
\end{equation}
The unsigned degenerate Stirling numbers of the first kind are given by
\begin{displaymath}
{n \brack k}_{\lambda}=(-1)^{n-k}S_{1,\lambda}(n,k),\quad (n, k\ge 0).
\end{displaymath}
Thus, by \eqref{12}, we get
\begin{equation}\label{13}
\frac{1}{k!}\big(-\log_{\lambda}(1-t)\big)^{k}=\frac{1}{k!} \log_{-\lambda}^k \left(\frac1{1-t}\right) =\sum_{n=k}^{\infty}
\begin{bmatrix}n\\k\end{bmatrix}_\lambda\ \frac{t^n}{n!},\quad (k\ge 0).
\end{equation} \par
In combinatorics, the term derangement refers to a permutation of a set of objects where none of the elements end up in their original position. A classical example is the hat-check problem: Imagine $n$ guests check their hats at a party, and at the end, the hat are returned randomly. A derangement is a scenario where no one gets their own hat back. The number of derangements for a set of $n$ elements is denoted by $d_{n}$ and called the $n$-th derangement number, also known as the subfactorial of $n$ or the $n$-th de Montmort number.
The $n$-th derangement number $d_n$ is given by
$$d_n =n! \sum_{k=0}^n\frac{(-1)^k}{k!},\quad(n\ge 0).$$
Note that
\begin{equation}\label{14}
\frac1{1-t}e^{-t} =\sum_{n=0}^\infty d_n\ \frac{t^n}{n!},\quad(\rm see\ [15]).
\end{equation}
Recently, Kim-Kim introduced the degenerate derangement numbers given by
\begin{equation}\label{15}
\frac1{1-t}e_\lambda^{-1}(t) =\sum_{n=0}^\infty d_{n,\lambda}\ \frac{t^n}{n!},\quad(\rm see\ [15]).
\end{equation}
Note that
  $$\lim_{\lambda\to 0}d_{n,\lambda}=d_n,\quad(n\ge 0).$$ \par
The beta function is defined by (see [1,12])
\begin{equation*}
B(z_1,z_2)=\int_{0}^{1}t^{z_1-1}(1-t)^{z_2-1}dt, \quad (\mathrm{Re}(z_1)>0, \ \mathrm{Re}(z_2)>0).
\end{equation*}
The beta function is related to the gamma function given by
\begin{equation}
B(z_1,z_2)=\frac{\Gamma(z_1)\Gamma(z_2)}{\Gamma(z_1+z_2)}. \label{15-1}
\end{equation}

\section{Degenerate harmonic numbers and some related numbers}

From \eqref{4} and \eqref{6}, we note that
\begin{equation}\label{16}
\begin{split}
\sum_{n=1}^\infty H_{n,\lambda} \ t^n
&= \frac{1}{1-t}\log_{-\lambda}\Big(\frac1{1-t}\Big) = \frac1{1-t}\log_{-\lambda}\left(1+ \frac{t}{1-t}\right)\\
& = \frac{1}{1-t}\sum_{k=1}^\infty \binom{-\lambda-1}{k-1} \frac{1}{k} \Big(\frac{t}{1-t}\Big)^{k} \\
& =\sum_{k=1}^{\infty}\frac{(-1)^{k-1}}{k}\binom{\lambda+k-1}{k-1}\sum_{n=0}^{\infty}\binom{n+k}{k}t^{n+k} \\
& =\sum_{k=1}^{\infty}\frac{(-1)^{k-1}}{k}\binom{\lambda+k-1}{k-1}\sum_{n=k}^{\infty}\binom{n}{k}t^{n} \\
& =\sum_{n=1}^\infty \sum_{k=1}^n \binom{n}{k} \frac{(-1)^{k-1}}{k} \binom{\lambda+k-1}{k-1} t^n.
\end{split}
\end{equation}
Therefore, by \eqref{16}, we obtain the following theorem.\par
\vskip 1cm
\begin{theorem}
For $n \in \mathbb{N}$, we have
\begin{equation}\label{17}
H_{0,\lambda}=0,\quad H_{n,\lambda}= \sum_{k=1}^n \binom{n}{k} \frac{(-1)^{k-1}}{k} \binom{\lambda+k-1}{k-1}.
\end{equation}
\end{theorem}

From \eqref{5}, we have
\begin{equation}\label{18}
H_{n,\lambda}-H_{n-1,\lambda}=\frac{(-1)^{n-1}}{n}\binom{\lambda-1}{n-1}, \quad (n \in \mathbb{N}),
\end{equation}
and, from \eqref{17}, we also have
\begin{equation}\label{19}
\begin{split}
H_{n,\lambda}-H_{n-1,\lambda}
&= \sum_{k=1}^n \binom{n}{k} \frac{(-1)^{k-1}}{k} \binom{\lambda+k-1}{k-1}-\sum_{k=1}^{n-1} \binom{n-1}{k} \frac{(-1)^{k-1}}{k} \binom{\lambda+k-1}{k-1} \\
&= \frac{(-1)^{n-1}}{n} \binom{\lambda+n-1}{n-1}+\sum_{k=1}^{n-1}\binom{n-1}{k-1}\frac{(-1)^{k-1}}{k}\binom{\lambda+k-1}{k-1}.
\end{split}
\end{equation}
Therefore, by combining \eqref{18} and \eqref{19}, we obtain the following corollary.
\par
\vskip 1cm
\begin{corollary}
For $n \in \mathbb{N}$, we have
\begin{equation*}
\frac{(-1)^{n-1}}{n!} \Big\{ (\lambda-1)_{n-1}- \langle \lambda+1 \rangle_{n-1}\Big\}=\sum_{k=1}^{n-1}\binom{n-1}{k-1}\frac{(-1)^{k-1}}{k}\binom{\lambda+k-1}{k-1},
\end{equation*}
where the rising factorial sequence is given by
\begin{equation*}
\langle x \rangle_{0}=1, \quad \langle x \rangle_{n}=x(x+1)\cdots (x+n-1),\quad (n \ge 1).
\end{equation*}
\end{corollary}
\vskip 0.5cm
The following binomial inversion theorem is well known (see [6,8,17,18]).
\begin{theorem}
The following relation holds true.
\begin{equation}\label{20}
a_n = \sum_{k=1}^n \binom{n}{k} (-1)^{k-1} b_k,\ (n \ge 1)\ \iff \ \ b_n = \sum_{k=1}^n \binom{n}{k} (-1)^{k-1} a_k, \ (n \ge 1).
\end{equation}
\end{theorem}
Assume that $a_n = \sum_{k=1}^n \binom{n}{k} (-1)^{k-1} b_k$ holds for all $n \in \mathbb{N}$. Then, by using \eqref{15-1}, we have
\begin{equation}\label{21}
\begin{split}
\sum_{k=1}^n \binom{n}{k}(-1)^{k-1}\frac{a_k}{k}
&=\sum_{k=1}^n \binom{n}{k} \frac{(-1)^{k-1}}{k} \sum_{l=1}^k\binom{k}{l}(-1)^{l-1}b_l\\
&= \sum_{l=1}^n (-1)^{l-1} b_l \sum_{k=l}^n \binom{n}{k} \frac{(-1)^{k-1}}{k} \binom{k}{l} \\
&= \sum_{l=1}^n (-1)^{l-1} b_l \binom{n}{l} \sum_{k=l}^n  \frac{(-1)^{k-1}}{k} \binom{n-l}{k-l} \\
&= \sum_{l=1}^n (-1)^{l-1} b_l \binom{n}{l} \sum_{k=0}^{n-l}  \frac{(-1)^{k+l-1}}{k+l} \binom{n-l}{k} \\
&= \sum_{l=1}^n b_l\binom{n}{l} \sum_{k=l}^{n-l} \binom{n-l}{k} \frac{(-1)^k}{k+l}\\
&= \sum_{l=1}^n b_l\binom{n}{l} \int_0^1 (1-t)^{n-l} t^{l-1}dt \\
&= \sum_{l=1}^n b_l\binom{n}{l}B(n-l+1,l) \\
&= \sum_{l=1}^n b_l\binom{n}{l}\frac{\Gamma(n-l+1)\Gamma(l)}{\Gamma(n+1)}\\
&=\sum_{l=1}^n \frac{b_l}{l}.
\end{split}
\end{equation}
Therefore, by combining \eqref{20} and \eqref{21}, we obtain the following lemma.
\par\vskip 1cm
\begin{lemma}
For any sequences $\{ a_{n}\}_{n=1}^{\infty}$ and $\{b_{n}\}_{n=1}^{\infty}$, the following holds true.
\begin{equation}\label{22}
\sum_{k=1}^n \binom{n}{k} (-1)^{k-1} a_k=b_n, \ (n \ge 1) \Rightarrow \sum_{k=1}^n \binom{n}{k}(-1)^{k-1}\frac{a_k}{k}=\sum_{k=1}^n \frac{b_k}{k}, \ (n \ge 1).
\end{equation}
\end{lemma}
\vskip 0.5cm
Applying \eqref{22} with $a_{k}$ replaced by $a_{k}^{\prime}=\frac{a_{k}}{k}$ and
$b_{n}$ replaced by $b_{n}^{\prime}=\sum_{k=1}^{n}\frac{b_{k}}{k}$, we obtain
\begin{align*}
\sum_{k=1}^{n}\binom{n}{k}(-1)^{k-1}\frac{a_{k}}{k^{2}}&=\sum_{k=1}^{n}\binom{n}{k}(-1)^{k-1}\frac{a_{k}^{\prime}}{k}
=\sum_{k_{2}=1}^{n}\frac{b_{k_{2}}^{\prime}}{k_{2}} \\
&=\sum_{k_{2}=1}^{n}\frac{1}{k_{2}}\sum_{k_{1}=1}^{k_{2}}\frac{b_{k_{1}}}{k_{1}}=\sum_{k_{2}=1}^{n}\sum_{k_{1}=1}^{k_{2}}\frac{b_{k_{1}}}{k_{1}k_{2}}. \nonumber
\end{align*}
Applying \eqref{22} with $a_{k}$ replaced by $a_{k}^{\prime \prime}=\frac{a_{k}}{k^{2}}$ and
$b_{n}$ replaced by $b_{n}^{\prime \prime}=\sum_{k_{2}=1}^{n}\sum_{k_{1}=1}^{k_{2}}\frac{b_{k_{1}}}{k_{1}k_{2}}$, we obtain
\begin{align*}
\sum_{k=1}^{n}\binom{n}{k}(-1)^{k-1}\frac{a_{k}}{k^{3}}&=\sum_{k=1}^{n}\binom{n}{k}(-1)^{k-1}\frac{a_{k}^{\prime \prime}}{k}
=\sum_{k_{3}=1}^{n}\frac{b_{k_{3}}^{\prime \prime}}{k_{3}} \\
&=\sum_{k_{3}=1}^{n}\frac{1}{k_{3}}\sum_{k_{2}=1}^{k_{3}}\sum_{k_{1}=1}^{k_{2}}\frac{b_{k_{1}}}{k_{1}k_{2}}=\sum_{k_{3}=1}^{n}\sum_{k_{2}=1}^{k_{3}}\sum_{k_{1}=1}^{k_{2}}\frac{b_{k_{1}}}{k_{1}k_{2}k_{3}}.
\nonumber
\end{align*}
Continuing this process, we have
\begin{align}
\sum_{k=1}^{n}\binom{n}{k}(-1)^{k-1}\frac{a_{k}}{k^{m}}&=\sum_{k_{m}=1}^{n}\sum_{k_{m-1}=1}^{k_{m}}\sum_{k_{m-2}=1}^{k_{m-1}} \cdots \sum_{k_{1}=1}^{k_{2}}\frac{b_{k_{1}}}{k_{1}k_{2}\cdots k_{m}} \label{23} \\
&=\sum_{1\le k_1 \le k_2\le \cdots \le k_m\le n} \frac{b_{k_{1}}}{k_{1}k_{2}\cdots k_{m}}. \nonumber
\end{align}
Therefore we obtain the following theorem from \eqref{23}.
\begin{theorem}
Let $\{ a_{n}\}_{n=1}^{\infty}$ and $\{b_{n}\}_{n=1}^{\infty}$ be any sequences. Assume that the following holds:
\begin{equation*}
\sum_{k=1}^n \binom{n}{k} (-1)^{k-1} a_k=b_n, \ (n \ge 1).
\end{equation*}
Then we have
\begin{equation*}
\sum_{k=1}^{n}\binom{n}{k}(-1)^{k-1}\frac{a_{k}}{k^{m}}=\sum_{1\le k_1 \le k_2\le \cdots \le k_m\le n} \frac{b_{k_{1}}}{k_{1}k_{2}\cdots k_{m}}, \ (n, \ m \ge 1).
\end{equation*}
\end{theorem}
Recall from \eqref{17} that $ \sum_{k=1}^n \binom{n}{k} (-1)^{k-1} \frac{\binom{\lambda+k-1}{k-1}}{k}=H_{n,\lambda}$. So we may apply Theorem 2.5 with $a_{k}=\frac{1}{k}\binom{\lambda+k-1}{k-1}$ and $b_{n}=H_{n,\lambda}$, and, by using \eqref{5}, get
\begin{align}
\sum_{k=1}^n \binom{n}{k}\frac{(-1)^{k-1}}{k^m}\binom{\lambda+k-1}{k-1}&=\sum_{k=1}^{n}\binom{n}{k}(-1)^{k-1}\frac{a_{k}}{k^{m-1}}
=\sum_{1 \le k_2 \le \cdots k_m \le n} \frac{H_{k_2,\lambda}}{k_2 \cdots k_m} \label{24} \\
&=\sum_{1 \le k_2 \le \cdots k_m \le n}\frac{1}{k_2 \cdots k_m}\sum_{k_1=1}^{k_2}\binom{\lambda -1}{k_1 -1}\frac{(-1)^{k_1-1}}{k_1} \nonumber  \\
&= \sum_{1\le k_1 \le k_2\le \cdots\le k_m\le n}    \frac{(-1)^{k_1 -1}}{k_1 k_2 \cdots k_m} \binom{\lambda-1}{k_1 -1}. \nonumber
\end{align}
Therefore, by \eqref{24}, we obtain the following theorem.
\par\vskip 1cm
\begin{theorem}
For $n,m\in \mathbb{N}$, we have
\begin{equation*}
\sum_{k=1}^n \binom{n}{k}\frac{(-1)^{k-1}}{k^m}\binom{\lambda+k-1}{k-1}
= \sum_{1\le k_1 \le k_2\le \cdots\le k_m\le n}    \frac{(-1)^{k_1 -1}}{k_1 k_2 \cdots k_m} \binom{\lambda-1}{k_1 -1}.
\end{equation*}
\end{theorem}
By \eqref{10}, we easily get
\begin{equation}\label{27}
\mathrm{Li}_{k,\lambda}(t)= \sum_{n=1}^\infty \frac{(-\lambda)^{n-1}\ (1)_{n,1/\lambda}}{(n-1)! \ n^k }\  t^n
=  \sum_{n=1}^\infty \frac{(-1)^{n-1}}{n^k} \binom{\lambda-1}{n-1}\  t^n .
\end{equation}
Define the numbers $K_{n,\lambda}^{(m)},\ (m, n\ge 1)$, by
\begin{equation}\label{28}
\frac{-1}{1-t}\mathrm{Li}_{m,-\lambda}\left(\frac{-t}{1-t}\right)=\sum_{n=1}^\infty K_{n,\lambda}^{(m)}\ t^n.
\end{equation}
We note here that $K_{n,\lambda}^{(1)}=H_{n,\lambda}$ (see \eqref{10-1}, \eqref{28}).
By \eqref{27}, we get
\begin{equation}\label{29}
\begin{split}
\mathrm{Li}_{m,-\lambda}\left(\frac{-t}{1-t}\right)&=\sum_{k=1}^\infty  \frac{(-1)^{k-1}}{k^m} \binom{-\lambda-1}{k-1}\left(\frac{-t}{1-t}\right)^k\\
&=-\sum_{k=1}^\infty  \frac{t^{k}}{k^m} \binom{-\lambda-1}{k-1}\sum_{n=0}^\infty \binom{n+k-1}{k-1}t^n\\
&=-\sum_{k=1}^\infty  \frac{t^{k}}{k^m} \binom{-\lambda-1}{k-1}\sum_{n=k}^\infty \binom{n-1}{k-1}t^{n-k}\\
&=-\sum_{n=1}^\infty \sum_{k=1}^n \binom{n-1}{k-1}\frac{1}{k^m} \binom{-\lambda-1}{k-1} t^{n}\\
&=-\sum_{n=1}^\infty \frac1n \sum_{k=1}^n \binom{n}{k}\frac{1}{k^{m-1}} \binom{-\lambda-1}{k-1}t^n.
\end{split}
\end{equation}
From \eqref{28} and \eqref{29}, we note that
\begin{equation}\label{30}
\begin{split}
\sum_{n=1}^\infty K_{n,\lambda}^{(m)}\ t^n
&= \frac{-1}{1-t}\mathrm{Li}_{m,-\lambda}\left(\frac{-t}{1-t}\right)\\
&=\sum_{l=0}^{\infty}t^{l}\sum_{j=1}^\infty\frac1j\sum_{k=1}^j \binom{j}{k}\frac{1}{k^{m-1}}\binom{-\lambda-1}{k-1}t^j\\
&=\sum_{n=1}^\infty \sum_{j=1}^n\frac1j \sum_{k=1}^j
 \binom{j}{k}\frac{(-1)^{k-1}}{k^{m-1}}\binom{\lambda +k-1}{k-1}t^n.
 \end{split}
\end{equation}
Thus, by \eqref{30} and Theorem 2.6, we get
\begin{equation}\label{31}
\begin{split}
 K_{n,\lambda}^{(m)}
&= \sum_{j=1}^n\frac1j \sum_{k=1}^j  \binom{j}{k}\frac{(-1)^{k-1}}{k^{m-1}}\binom{\lambda+k-1}{k-1}\\
&= \sum_{k_m=1}^n\frac{1}{k_m} \sum_{k=1}^{k_m} \binom{k_m}{k}\frac{(-1)^{k-1}}{k^{m-1}}\binom{\lambda+k-1}{k-1}\\
&= \sum_{k_m=1}^n\frac{1}{k_m}\sum_{1\le k_1 \le k_2\le \cdots\le k_{m-1}\le k_m}    \frac{(-1)^{k_1 -1}}{k_1 k_2 \cdots k_{m-1}} \binom{\lambda-1}{k_1 -1}\\
&= \sum_{1\le k_1\le k_2 \le \cdots \le k_{m} \le n}\frac{(-1)^{k_1-1}}{k_1 k_2 \cdots k_{m}}    \binom{\lambda-1}{k_1-1}.
\end{split}
\end{equation}
Therefore, by \eqref{31}, we obtain the following theorem.
\par
\vskip 1cm
\begin{theorem}
For $n,m\in \mathbb{N}$, we have
\begin{equation*}
 K_{n,\lambda}^{(m)}
 = \sum_{1\le k_1\le k_2 \le \cdots \le k_m \le n}\frac{(-1)^{k_1-1}}{k_1 k_2 \cdots k_{m}}    \binom{\lambda-1}{k_1-1}.
\end{equation*}
 \end{theorem}

 \vskip 0.5cm
 By \eqref{8} and \eqref{29}, we get
 \begin{equation}\label{32}
 \begin{split}
\frac{-1}{1-t}\mathrm{Li}_{m,-\lambda}\left(\frac{-t}{1-t}\right)
&=\frac1{1-t}\sum_{k=1}^\infty \frac{1}{k^m} \binom{-\lambda-1}{k-1}{k!}\frac1{k!}\left(\frac{t}{1-t}\right)^k\\
&=\frac1{1-t}\sum_{k=1}^\infty \frac{(-1)^{k-1}}{k^m} \binom{\lambda+k-1}{k-1}{k!}\sum_{j=k}^\infty L(j,k)\frac{t^j}{j!}\\
&=\sum_{l=0}^\infty t^l \sum_{j=1}^\infty \left( \frac1{j!} \sum_{k=1}^j
\frac{(-1)^{k-1}}{k^m} \binom{\lambda+k-1}{k-1}{k!}\ L(j,k)\right)t^j \\
&=\sum_{n=1}^\infty \left(\sum_{j=1}^n \frac1{j!} \sum_{k=1}^j
\frac{(-1)^{k-1}}{k^m} \binom{\lambda+k-1}{k-1}{k!}\ L(j,k)\right)t^n.
\end{split}
\end{equation}
From \eqref{28} and \eqref{32}, we obtain the following theorem.
\par
\vskip 1cm
\begin{theorem}
For $m,n \in \mathbb{N}$, we have
\begin{equation*}
K_{n,\lambda}^{(m)}=\sum_{j=1}^n \frac1{j!} \sum_{k=1}^j
\frac{(-1)^{k-1}}{k^m} \binom{\lambda+k-1}{k-1}{k!}\ L(j,k).\\
\end{equation*}
\end{theorem}

\vskip 0.5cm
Now, we define the degenerate harmonic numbers of order $k$ by
\begin{equation}\label{33}
\frac1{1-t}\mathrm{Li}_{k,\lambda}(t) = \sum_{n=1}^\infty H_{n,\lambda}^{(k)}\ t^n, \quad H_{0,\lambda}^{(k)}=0.
\end{equation}
When $k=1$, we have $ H_{n,\lambda}^{(k)}=H_{n,\lambda}$ (see \eqref{10-1}, \eqref{33}). From \eqref{27} and \eqref{33}, we note that
\begin{equation}\label{34}
\begin{split}
\sum_{n=1}^\infty H_{n,\lambda}^{(k)}\ t^n &=\frac1{1-t}\mathrm{Li}_{k,\lambda}(t)= \sum_{l=0}^\infty t^l \sum_{m=1}^{\infty}\frac{(-1)^{m-1}}{m^{k}}\binom{\lambda-1}{m-1} t^m \\
&= \sum_{n=1}^\infty \sum_{m=1}^n \frac{(-1)^{m-1}}{m^k}\binom{\lambda-1}{m-1} t^n.
\end{split}
\end{equation}
Therefore, by \eqref{34}, we obtain the following theorem.
\par
\vskip 1cm
\begin{theorem}
For $n,k\in \mathbb{N}$, we have
\begin{equation*}
H_{n,\lambda}^{(k)}= \sum_{m=1}^n  \frac{(-1)^{m-1}}{m^k} \binom{\lambda-1}{m-1}.
\end{equation*}
\end{theorem}
Here we note that
$$\lim_{\lambda\to 0}H_{n,\lambda}^{(k)}=\sum_{m=1}^n\frac{1}{m^k}=H_n^{(k)}.$$
By \eqref{6}, we get
\begin{equation}\label{35}
\begin{split}
\sum_{n=1}^\infty H_{n,\lambda} t^n
 &= \frac1{1-t}\log_{-\lambda}\left(\frac{1}{1-t}\right)\\
&=e_{-\lambda}\left(\log_{-\lambda}\left(\frac1{1-t}\right)\right) \log_{-\lambda}\left(\frac1{1-t}\right)\\
&= \sum_{k=0}^\infty \frac{(1)_{k,-\lambda}}{k!} \log_{-\lambda}^k\left(\frac1{1-t}\right) \log_{-\lambda}\left(\frac1{1-t}\right)\\
&= \sum_{k=0}^\infty \frac{(1)_{k,-\lambda}}{k!}(k+1)! \frac1{(k+1)!}   \log_{-\lambda}^{k+1}\left(\frac1{1-t}\right)\\
&= \sum_{k=0}^\infty (k+1)  (1)_{k,-\lambda} \sum_{n=k+1}^\infty    \begin{bmatrix}n\\k+1\end{bmatrix}_\lambda \frac{t^n}{n!}\\
&= \sum_{k=1}^\infty k \ (1)_{k-1,-\lambda} \sum_{n=k}^\infty    \begin{bmatrix}n\\k\end{bmatrix}_\lambda \frac{t^n}{n!}\\
&= \sum_{n=1}^\infty \sum_{k=1}^n k \  (1)_{k-1,-\lambda}  \begin{bmatrix}n\\k\end{bmatrix}_\lambda \frac{t^n}{n!}.
\end{split}
\end{equation}
Therefore, by \eqref{35}, we obtain the following theorem.
\par
\vskip 1cm
\begin{theorem}
For $n \in \mathbb{N}$, we have
\begin{equation*}
H_{n,\lambda}= \sum_{k=1}^n  k\  (1)_{k,-\lambda}\begin{bmatrix}n\\ k\end{bmatrix}_\lambda.
\end{equation*}
 \end{theorem}
\begin{remark}
By \eqref{27} and Theorem 2.9, we obtain
\begin{equation*}
\mathrm{Li}_{k,\lambda}(t)=\sum_{n=1}^{\infty}\frac{(-1)^{n-1}}{n^{k}}\binom{\lambda-1}{n-1}t^{n}
=\sum_{n=1}^{\infty}\big(H_{n,\lambda}^{(k)}-H_{n-1,\lambda}^{(k)}\big)t^{n}.
\end{equation*}
In particular, for $k=1$, we have
\begin{align}
\log_{-\lambda}\Big(\frac{1}{1-t}\Big)&=Li_{1,\lambda}(t)=\sum_{n=1}^{\infty}\frac{(-1)^{n-1}}{n}\binom{\lambda-1}{n-1}t^{n} \label{36} \\
&=\sum_{n=1}^{\infty}\big(H_{n,\lambda}-H_{n-1,\lambda}\big)t^{n}. \nonumber
\end{align}
\end{remark}

\vskip 0.5cm
Now, from \eqref{6}, \eqref{15} and \eqref{36}, we observe that
\begin{equation}\label{37}
\begin{split}
\sum_{n=1}^\infty H_{n,\lambda}\  t^n
&= \frac1{1-t}\log_{-\lambda}\left(\frac{1}{1-t}\right)
= \frac1{1-t} e_\lambda^{-1}(t) \log_{-\lambda}\left(\frac{1}{1-t}\right) e_\lambda(t)\\
&= \sum_{l=0}^\infty d_{l,\lambda}\frac{t^l}{l!}
\sum_{j=1}^\infty \left(\sum_{k=1}^j \left(H_{k,\lambda}-H_{k-1,\lambda}\right)\ (1)_{j-k,\lambda}\frac1{(j-k)!}\right)t^j\\
&=  \sum_{n=1}^\infty \sum_{j=1}^n\binom{n}{j}d_{n-j,\lambda}
\sum_{k=1}^j k! \binom{j}{k} \left(H_{k,\lambda}-H_{k-1,\lambda}\right) \ (1)_{j-k,\lambda}\frac{t^n}{n!}.
\end{split}
\end{equation}
Therefore, by \eqref{37}, we obtain the following theorem.
\par
\vskip 1cm
\begin{theorem}
For $n \in \mathbb{N}$, we have
\begin{equation*}
  H_{n,\lambda}=  \frac{1}{n!} \sum_{j=1}^n\binom{n}{j}d_{n-j,\lambda}
\sum_{k=1}^j k! \binom{j}{k} \left(H_{k,\lambda}-H_{k-1,\lambda}\right)\ (1)_{j-k,\lambda}.
\end{equation*}
\end{theorem}

\section{Conclusion}
In this paper, we derived some expressions for the degenerate harmonic numbers $H_{n,\lambda}$, and for the degenerate harmonic numbers of order $m$, $H_{n,\lambda}^{(m)}$. Indeed, we found the following:
\begin{align*}
H_{n,\lambda}&= \sum_{k=1}^n  \frac{(-1)^{k-1}}{k} \binom{\lambda-1}{k-1} \\
&=\sum_{k=1}^n \binom{n}{k} \frac{(-1)^{k-1}}{k} \binom{\lambda+k-1}{k-1} \\
&= \sum_{k=1}^n  k\  (1)_{k,-\lambda}\begin{bmatrix}n\\ k\end{bmatrix}_\lambda \\
&= \frac{1}{n!} \sum_{j=1}^n\binom{n}{j}d_{n-j,\lambda}
\sum_{k=1}^j k! \binom{j}{k} \left(H_{k,\lambda}-H_{k-1,\lambda}\right)\ (1)_{j-k,\lambda},
\end{align*}
\begin{equation*}
H_{n,\lambda}^{(m)}= \sum_{k=1}^n  \frac{(-1)^{k-1}}{k^m} \binom{\lambda-1}{k-1}.
\end{equation*} \par
In addition, we introduced the numbers $K_{n,\lambda}^{(m)}$ and showed some explicit expressions for them, which are closely connected with $H_{n,\lambda}^{(m)}$ and reduce to $K_{n,\lambda}^{(1)}=H_{n,\lambda}$, for $m=1$. Indeed, we have deduced
\begin{align*}
K_{n,\lambda}^{(m)}&=\sum_{k=1}^n \binom{n}{k}\frac{(-1)^{k-1}}{k^m}\binom{\lambda+k-1}{k-1} \\
&=\sum_{j=1}^n \frac1{j!} \sum_{k=1}^j
\frac{(-1)^{k-1}}{k^m} \binom{\lambda+k-1}{k-1}{k!}\ L(j,k) \\
&= \sum_{1\le k_1 \le k_2\le \cdots\le k_m\le n}\frac{(-1)^{k_1 -1}}{k_1 k_2 \cdots k_m} \binom{\lambda-1}{k_1 -1}.
\end{align*} \par

\end{document}